# PERFECT SIMULATION FOR A CLASS OF POSITIVE RECURRENT MARKOV CHAINS

By Stephen B. Connor and Wilfrid S. Kendall

*University of Warwick*

This paper generalizes the work of Kendall [*Electron. Comm. Probab.* **9** (2004) 140–151], which showed that perfect simulation, in the form of dominated coupling from the past, is always possible (although not necessarily practical) for geometrically ergodic Markov chains. Here, we consider the more general situation of positive recurrent chains and explore when it is possible to produce such a simulation algorithm for these chains. We introduce a class of chains which we name *tame*, for which we show that perfect simulation is possible.

**1. Introduction.** Perfect simulation was first introduced by Propp and Wilson [20] as a method for sampling from the exact stationary distribution of an ergodic Markov chain. Foss and Tweedie [7] showed that this classic coupling from the past (*CFTP*) algorithm is possible (in principle, if not in practice) if and only if the Markov chain is uniformly ergodic.

More recently, Kendall [14] showed that all geometrically ergodic chains possess (again, possibly impractical) *dominated CFTP* algorithms (as introduced in [13, 16]). This suggests the questions: what if $X$ is *subgeometrically ergodic*? Might it be the case that *all* positive recurrent Markov chains possess (impractical) *domCFTP* algorithms?

In this paper, we introduce a new class of positive-recurrent chains (*tame chains*) for which *domCFTP* is possible in principle.

Note that the *practicality* of *CFTP* algorithms is subject to a number of interesting constraints: methods using coadapted coupling will deliver answers at a slower exponential rate than ordinary Markov chain Monte Carlo for many chains [1, 19]; in general, the coalescence of paths from many different starting states (an intrinsic feature of *CFTP*) may be expected to be slower than pairwise coupling; finally, the theory of randomized algorithms can be used to demonstrate the existence of problems for which there will not even be any fully-polynomial randomized approximation schemes (subject to the complexity theory assumption $RP \neq NP$; Jerrum [12] discusses results of this nature for counting algorithms for independent sets).

Considerations of the practicality of *CFTP* raise many further interesting research questions; however, in this paper, we focus on considering whether (for all Markov chains with a specified property) there can exist *domCFTP* algorithms, practical or not.

To make this a meaningful exercise, it is necessary to be clearer about what one is allowed to do as part of an impractical algorithm. The [7] result for uniform ergodicity presumes that one is able to identify when regeneration occurs for the target Markov chain subsampled every $k$ time steps (where $k$ is the order of the whole state space considered as a small set of the chain) and that one can then draw from the regeneration distribution and the $k$-step transition probability distribution conditioned on no regeneration. One must assume more in order to cover the geometrically ergodic case [14], namely, that it is possible to couple the target chain and the dominating chain when subsampled every $k$ time steps, preserving the domination while so doing. Here, $k$ is the order of a small set for a particular Foster–Lyapunov criterion for the geometric ergodicity property. In fact, something more must also be assumed: it must be possible to implement the coupling between target chain and dominating process in a monotonic fashion, even when conditioning on small-set regeneration occurring or not occurring. In fact, we do not need to assume any more than this when dealing with the tame chains introduced below, except that the subsampling order $k$ is now not fixed for all time, but can vary according to the current value of the dominating process.

The impracticality of these *CFTP* algorithms thus has two aspects. First, the question of expected run time is not addressed at all. Second, for the most part, the assumptions described above amount to supposing that we can translate into practice the theoretical possibility of implementing various stochastic dominations as couplings (guaranteed by theory expounded in, e.g., [17], Chapter IV). However, it should be noted that practical and implemented *CFTP* algorithms can correspond very closely to these general schemes. For example, the *CFTP* algorithm resulting from the result of Foss and Tweedie [7] is essentially the simplest case of the exact sampling algorithm proposed by Green and Murdoch [9]; the scheme proposed in [14] is closely related to fast *domCFTP* algorithms for perpetuities with sample step $k = 1$.

In this paper, we investigate the problems that occur in the move from geometric to subgeometric ergodicity. We begin by recalling some useful results concerning rates of ergodicity. Section 2 then reviews the result of [14].

The bulk of the new material in this paper is to be found in Section 3. There, we introduce the notion of a *tame* chain (Definition 14) and demonstrate that *domCFTP* is possible for such chains (Theorem 15). A description of the *domCFTP* algorithm for tame chains is provided in Section 3.3; the reader is referred to [15] for an introduction to the classical form of *domCFTP*. We also prove some sufficient conditions for a polynomially ergodic chain to be tame (Theorems 21 and 22). However, these conditions are not necessary; Section 4.4 contains an example of a polynomially ergodic chain which does not satisfy these conditions and yet is still tame. The existence of a polynomially ergodic chain that is not tame is currently an open question.

1.1. *Definitions and notation.* Let $X = (X_0, X_1, \ldots)$ be a discrete-time Markov chain on a Polish state space $\mathcal{X}$. The Markov transition kernel for $X$ is denoted by $P$ and the $n$-step kernel by $P^n$:

$$P^n(x, E) = \mathbb{P}_x[X_n \in E],$$

where $\mathbb{P}_x$ is the conditional distribution of the chain given $X_0 = x$. The corresponding expectation operator will be denoted by $\mathbb{E}_x$. If $g$ is a nonnegative function, then we write $Pg(x)$ for the function $\int g(y) P(x, dy)$ and for a signed measure $\mu$, we write $\mu(g)$ for $\int g(y) \mu(dy)$. The $f$-norm is defined as $\|\mu\|_f := \sup_{g : |g| \leq f} |\mu(g)|$; taking $f \equiv 1$ yields the usual total variation norm, for which we will simply write $\|\mu\|$.

We assume throughout that $X$ is aperiodic (in the sense of [18]) and Harris-recurrent. The stationary distribution of $X$ shall be denoted by $\pi$ and the *first hitting time* of a measurable set $A \subseteq \mathcal{X}$ by $\tau_A = \min\{n \geq 1 : X_n \in A\}$.

The notion of *small sets* will feature heavily throughout this paper.

DEFINITION 1. A subset $C \subseteq \mathcal{X}$ is a small set (of order $m$) for the Markov chain $X$ if the following *minorization* condition holds: for some $\varepsilon \in (0, 1]$ and a probability measure $\nu$,

(1)     $\mathbb{P}_x[X_m \in E] \geq \varepsilon \nu(E)$     for all $x \in C$ and measurable $E \subset \mathcal{X}$.

In this case, we say that $C$ is *m-small*. Many results in the literature are couched in terms of the more general idea of *petite sets*; however, for aperiodic $\phi$-irreducible chains, the two notions are equivalent ([18], Theorem 5.5.7). Small sets allow the use of coupling constructions: specifically, if $X$ hits the small set $C$ at time $n$, then there is a positive chance ($\varepsilon$) that it regenerates at time $n + m$ (using the measure $\nu$). Furthermore, if regeneration occurs, then a single draw from $\nu$ may be used for any number of copies of $X$ belonging to $C$ at time $n$, resulting in their coalescence at time $n + m$. Small sets belong to a larger class of *pseudo-small* sets, as introduced in [22], but such sets only allow for the coupling of pairs of chains. Implementation of *domCFTP* requires a positive chance of a continuum of chains coalescing when belonging to a given set $C$, so we shall henceforth deal solely with small sets.

1.2. *Geometric ergodicity.* We first outline some relevant theory for geometrically ergodic chains.

DEFINITION 2. The chain $X$ is said to be *geometrically ergodic* if there exists a constant $\gamma \in (0,1)$ and some function $\Lambda : \mathcal{X} \to [0, \infty)$ such that for all $x$ in a full and absorbing set,

$$\|P^n(x, \cdot) - \pi(\cdot)\| \leq \Lambda(x) \gamma^n. \tag{2}$$

If $\Lambda$ can be chosen to be bounded, then $X$ is said to be *uniformly ergodic*.

Uniform ergodicity of $X$ can be shown to be equivalent to the whole state space $\mathcal{X}$ being a small set, in which case at every Markov chain step, there is a positive chance of coalescence, whereby chains started at all elements of the state space become equal simultaneously. Foss and Tweedie [7] use this to show that uniform ergodicity is equivalent to the existence of a *CFTP* algorithm for $X$ in the sense of Propp and Wilson [20].

The most common way to establish geometric ergodicity of a chain $X$ is to check the following geometric *Foster–Lyapunov* condition [8].

CONDITION GE. There exist positive constants $\beta < 1$ and $b < \infty$, a small set $C$ and a *scale function* $V : \mathcal{X} \to [1, \infty)$, bounded on $C$, such that

$$\mathbb{E}[V(X_{n+1})|X_n = x] \leq \beta V(x) + b\mathbf{1}_C(x). \tag{3}$$

Inequality (3) will be referred to as $\mathrm{GE}(V, \beta, b, C)$ when we need to be explicit about the scale function and constants. For simplicity, we will also often write inequality (3) as $PV \leq \beta V + b\mathbf{1}_C$. Under our global assumptions on $X$, this drift condition is actually equivalent to $X$ being geometrically ergodic ([18], Theorem 15.0.1). Furthermore, if $X$ satisfies (3), then we can take $\Lambda = V$ in equation (2).

Condition GE quantifies the way in which the chain $V(X)$ behaves as a supermartingale before $X$ hits $C$. When the chain hits $C$, it can then increase in expectation, but only by a bounded amount. The following result can be extracted from [18], Theorems 15.0.1 and 16.0.1.

THEOREM 3. *Suppose that $X$ is $\phi$-irreducible and aperiodic. Then $X$ is geometrically ergodic if and only if there exists $\kappa > 1$ such that the corresponding geometric moment of the first return time to $C$ is bounded, that is,*

$$\sup_{x \in C} \mathbb{E}_x[\kappa^{\tau_C}] < \infty. \tag{4}$$

The first hitting time of $C$ is related to drift conditions in the following way (extracted from [18], Theorem 11.3.5).

THEOREM 4. *For an ergodic chain $X$, the function $V_C(x) = \mathbb{E}_x[\tau_C]$ is the pointwise minimal solution to the inequality*

$$PV(x) \leq V(x) - 1, \qquad x \notin C. \tag{5}$$

Equation (5) is clearly a weaker drift condition than Condition GE and is equivalent to positive recurrence of $X$ [18]. It can be shown that (5) implies that all sublevel sets are small [18] and since $V$ is bounded on $C$, we will always take $C$ to be a sublevel set of the form $\{x \in \mathcal{X} : V(x) \leq d\}$.

We now present a couple of easy results concerning geometrically ergodic chains, which will prove to be of great importance later on. The first demonstrates how the scale function $V$ in (3) may be changed to obtain a new drift condition using the same small set.

LEMMA 5. *If the chain $X$ satisfies Condition $\mathrm{GE}(V, \beta, b, C)$, then for any $\xi \in (0, 1]$,*

$$PV^\xi \leq (\beta V)^\xi + b^\xi \mathbf{1}_C.$$

*Thus $\mathrm{GE}(V, \beta, b, C)$ implies $\mathrm{GE}(V^\xi, \beta^\xi, b^\xi, C)$.*

PROOF. Calculus shows that $(x+y)^\xi \leq x^\xi + y^\xi$ for $x, y \geq 0$ and $0 \leq \xi \leq 1$. The result follows by Jensen's inequality for $(PV)^\xi$, using (3). □

The second result shows that a geometric drift condition persists if we subsample the chain at some randomized stopping time.

LEMMA 6. *Suppose $X$ satisfies Condition $\mathrm{GE}(V, \beta, b, C)$. Then for any positive, integer-valued stopping time $\sigma$ (adapted to the natural filtration generated by $X$), we have*

$$\mathbb{E}_x[V(X_\sigma)] \leq \beta V(x) + b_1 \mathbf{1}_{C_1}(x),$$

*where $b_1 = b/(1-\beta)$ and $C_1 = \{x : V(x) \leq b/(\beta(1-\beta)^2)\} \cup C$.*

The same $\beta$, $b_1$ and $C_1$ work for all values of $\sigma$ since the constant $b_1$ absorbs the higher-order terms in $\beta$ below.

PROOF OF LEMMA 6. Iterate the drift condition (3) and treat the cases $\{\sigma = 1\}$ and $\{\sigma > 1\}$ separately:

$$\mathbb{E}_x[V(X_\sigma)] \leq \mathbb{E}_x\left[\beta^\sigma V(x) + b \sum_{j=1}^\sigma \beta^{j-1} \mathbf{1}_C(X_{\sigma-j})\right]$$

$$\leq (\beta V(x) + b\mathbf{1}_C(x))\mathbb{P}_x[\sigma = 1] + \left(\beta^2 V(x) + \frac{b}{1-\beta}\right)\mathbb{P}_x[\sigma > 1]$$

$$\leq (\beta V(x) + b\mathbf{1}_C(x))\mathbb{P}_x[\sigma = 1] + (\beta V(x) + b_1 \mathbf{1}_{C_1}(x))\mathbb{P}_x[\sigma > 1]$$

$$\leq \beta V(x) + b_1 \mathbf{1}_{C_1}(x). \qquad \square$$

1.3. *Polynomial ergodicity.* We now turn to polynomially ergodic chains and state some results which will prove useful in Section 3.4.

DEFINITION 7. The chain $X$ is said to be *polynomially ergodic* if there exists $\gamma \geq 0$ such that for all $x$ in a full and absorbing set,

$$(6) \qquad n^{\gamma} \|P^n(x, \cdot) - \pi(\cdot)\| \to 0 \qquad \text{as } n \to \infty.$$

As with geometric ergodicity, there is a Foster–Lyapunov drift condition that can be shown [11] to imply polynomial ergodicity (although the two are not equivalent in this case).

CONDITION PE. There exist constants $\alpha \in (0,1)$ and $b, c \in (0, \infty)$, a small set $C$ and a scale function $V : \mathcal{X} \to [1, \infty)$, bounded on $C$, such that

$$(7) \qquad \mathbb{E}[V(X_{n+1})|X_n = x] \leq V(x) - cV^{\alpha}(x) + b\mathbf{1}_C(x).$$

We will refer to (7) as PE$(V, c, \alpha, b, C)$ when we need to be explicit about the scale function and constants.

This drift condition again tells us that $V(X)$ behaves as a supermartingale before $X$ hits $C$, but that the drift toward the small set now occurs at a subgeometric rate (and hence $\tau_C$ has no exponential moment). Note that if $\alpha = 1$, then we regain Condition GE [for $c \in (0,1)$] and that we do not include the case $\alpha = 0$ here, for which the drift condition is equivalent to $X$ being simply positive recurrent.

Polynomially ergodic chains satisfy a result analogous to Lemma 5, with a similar proof ([11], Lemma 3.5).

LEMMA 8. *If the chain $X$ satisfies Condition* PE, *then for any $\xi \in (0, 1]$, there exists $0 < b_1 < \infty$ such that*

$$PV^{\xi} \leq V^{\xi} - c\xi V^{\alpha+\xi-1} + b_1 \mathbf{1}_C.$$

Note that as in Lemma 5, the same small set $C$ appears in the new drift condition when we change scale function in this way.

COROLLARY 9. *Suppose that $X$ satisfies Condition* PE. *Then for $x \notin C$,*

$$\mathbb{E}_x[\tau_C] \leq \frac{V^{1-\alpha}(x)}{c(1-\alpha)}.$$

PROOF. Set $\xi = 1 - \alpha$ in Lemma 8 to obtain

$$PV^{1-\alpha}(x) \leq V^{1-\alpha}(x) - c(1-\alpha) \qquad \text{for } x \notin C.$$

The result then follows from Theorem 4. □

Note, however, that there is no analogue to Lemma 6 (even if $\sigma$ is deterministic), since the geometric ergodicity case makes essential use of the convergence of the series $\sum \beta^j$.

The drift condition (7) can actually be shown to imply much more than the convergence in (6). From Theorem 3.6 of [11] we obtain the following which will be used in the proof of Theorem 22.

PROPOSITION 10. *Suppose $X$ satisfies Condition* PE. *Define, for each $1 \leq \rho \leq 1/(1-\alpha)$,*

$$V_\rho(x) = V^{1-\rho(1-\alpha)}(x) \quad \text{and} \quad r_\rho(n) = (n+1)^{\rho-1}. \tag{8}$$

*Then there exists a constant $M < \infty$ such that*

$$\mathbb{E}_x\left[\sum_{n=0}^{\tau_C - 1} r_\rho(n) V_\rho(X_n)\right] \leq MV(x). \tag{9}$$

Furthermore, from [4] we see that an upper bound for $M$ can be obtained directly from the drift condition (7).

**2. Geometric ergodicity implies domCFTP.** We now give a brief overview of the proof that all geometrically ergodic chains possess (not necessarily practical) *domCFTP* algorithms [14]. Recall that coadaptive coupling of Markov chains means that both chains have a common past expressed by a fixed filtration of $\sigma$-algebras.

DEFINITION 11. Suppose that $V$ is a scale function for a Harris-recurrent Markov chain $X$. We say that the stationary ergodic random process $Y$ on $[1, \infty)$ is a *dominating process for $X$ based on the scale function $V$* (with *threshold $h$* and *coalescence probability $\varepsilon$*) if it can be coupled coadaptively to realizations of $X^{x,-t}$ (the Markov chain $X$ begun at $x$ at time $-t$) as follows:

(a) for all $x \in \mathcal{X}$, $n > 0$ and $-t \leq 0$, almost surely

$$V(X^{x,-t}_{-t+n}) \leq Y_{-t+n} \implies V(X^{x,-t}_{-t+n+1}) \leq Y_{-t+n+1}; \tag{10}$$

(b) if $Y_n \leq h$, then the probability of *coalescence* at time $n+1$ is at least $\varepsilon$, where coalescence at time $n+1$ means that the set

$$\{X^{x,-t}_{n+1} : -t \leq n \text{ and } V(X^{x,-t}_n) \leq Y_n\} \tag{11}$$

is a singleton set;

(c) $\mathbb{P}[Y_n \leq h]$ is positive.

The following theorem is the main result of [14].

THEOREM 12. *If $X$ satisfies the drift condition*

$$PV \leq \beta V + b\mathbf{1}_C$$

*for $0 < \beta < 1$, then there exists a domCFTP algorithm for $X$ (possibly subject to subsampling) using a dominating process based on the scale $V$.*

The idea behind the proof of Theorem 12 is that a dominating process $Y$ satisfying equation (10) may be obtained by using Markov's inequality and the geometric drift condition for $X$. The result is that any chain satisfying Condition GE can be dominated by $Y = (d + b/\beta)\exp(U)$, where $U$ is the system workload of a $D/M/1$ queue, sampled at arrivals, with arrivals every $\log(1/\beta)$ units of time and service times being independent and of unit-rate exponential distribution. $U$ is positive recurrent only if $\beta < e^{-1}$, but a new geometric drift condition with $\beta$ replaced by $\beta^{k-1}$ can be produced by subsampling $X$ with a fixed subsampling period $k$; the proof uses the ideas of Lemma 6. If $k$ is chosen sufficiently large to fix $\beta^{k-1} < e^{-1}$, then the above argument produces a stationary dominating process for the subsampled chain.

Note that $Y$ is easy both to sample from in equilibrium and to run in reversed time, which is essential for implementation of *domCFTP*. Also, note that $Y$ belongs to a family of *universal* dominating processes for geometrically ergodic chains, although this dominator need not generally lead to a practical simulation algorithm. As noted in the introduction, the main difficulties in application are in implementing practical domination and in determining whether or not regeneration has occurred when $Y$ visits the set $\{Y \leq h\}$. This task is rendered even less practical if subsampling has taken place, since then, detailed knowledge of convolutions of the transition kernel for $X$ is required.

**3. domCFTP for suitable positive recurrent chains.** Theorem 12 leads to an obvious question: does there exist a similar *domCFTP* algorithm for chains not satisfying Condition GE? [Note that if we try to use the drift condition (7)—as in the proof of Theorem 12—to produce a dominating process for polynomially ergodic chains, then the resulting process is nonrecurrent.] In this section, we introduce a class of chains which possess a *domCFTP* algorithm.

The principal idea behind the subsequent work is to investigate when it is possible to subsample $X$ to produce a geometrically ergodic chain. For non-geometrically ergodic chains, a fixed subsampling interval will not work and so we seek an appropriate simple adaptive subsampling scheme. A similar scheme can then be used to *delay* the dominating process $Y$ constructed in Section 2 and to show that this new process $D$ dominates the chain $V(X)$ at the times when $D$ moves.

Several issues must be addressed in order to derive a *domCFTP* algorithm using this idea:

1. What is an appropriate adaptive subsampling scheme?
2. When does such a scheme exist?
3. How does the dominating process $D$ dominate $V(X)$ when $D$ moves?
4. Can we simulate $D$ in equilibrium and in reversed time?

The answers to these questions are quite subtle.

3.1. *Adaptive subsampling.* We begin by defining more carefully what we mean by an adaptive subsampling scheme.

DEFINITION 13. An *adaptive subsampling scheme* for the chain $X$ with respect to a scale function $V$ is a sequence of stopping times $\{\theta_n\}$ defined recursively by

(12) $$\theta_0 = 0; \qquad \theta_{n+1} = \theta_n + F(V(X_{\theta_n})),$$

where $F:[1,\infty) \to \{1,2,\ldots\}$ is a deterministic function.

Note that a set of stopping times $\{\theta_n\}$ such that $\{X_{\theta_n}\}$ is *uniformly ergodic* can be produced as follows. Using the Athreya–Nummelin split-chain construction [18], we may suppose that there is a state $\omega$ with $\pi(\omega) > 0$. Define

(13) $$F(V(x)) = \min\left\{m > 0 : \mathbb{P}_x[X_m = \omega] > \frac{\pi(\omega)}{2}\right\}.$$

Then the time until $\{X_{\theta_n}\}$ hits $\omega$ from any starting state $x$ is majorized by a geometric random variable with success probability $\pi(\omega)/2$. This implies that the subsampled chain is uniformly ergodic, as claimed. $F$ as defined in (13) depends upon knowledge of $\pi$, however, and we obviously do not have this available to us (it is the distribution from which we are trying to sample!). This example shows that adaptive subsampling can have drastic effects on $X$. However, construction of a *domCFTP* algorithm for $X$ using this subsampling scheme (in a manner to be described in Section 3.3) turns out to be impossible unless $X$ is itself uniformly ergodic.

Reverting to the previous discussion, suppose that there is an explicit adaptive subsampling scheme such that the chain $X' = \{X_{\theta_n}\}$ satisfies Condition GE with drift parameter $\beta < e^{-1}$. Then a candidate dominating process $D$ can be produced for $V(X)$ in the following way. Begin with an exponential queue workload process $Y$ that dominates $V(X')$ (as in Section 2). Then *slow down* $Y$ by generating *pauses* using some convenient function $S$ satisfying $S(z) \geq F(z')$ whenever $z \geq z'$, to produce the process $D$. That is, given $D_0 = Y_0 = z$, pause $D$ by setting

$$D_1 = D_2 = \cdots = D_{S(z)-1} = z.$$

Then define the law of $D_{S(z)}$ by $\mathcal{L}(D_{S(z)}|D_{S(z)-1} = z) = \mathcal{L}(Y_1|Y_0 = z)$. Iteration of this construction leads to a sequence of times $\{\sigma_n\}$ at which $D$ moves, defined recursively by

$$\sigma_{n+1} = \sigma_n + S(D_{\sigma_n}),$$

with $D$ constant on each interval of the form $[\sigma_n, \sigma_{n+1})$.

Such a process $D$ is a plausible candidate for a dominating process. To be suitable for use in a *domCFTP* algorithm, however, it must be possible to compute its equilibrium distribution. Now, $D$ as we have just defined it is only a semi-Markov process: it is Markovian at the times $\{\sigma_n\}$, but not during the delays between jumps. To remedy this, we augment the chain by adding a second coordinate $N$ that measures the time until the next jump of $D$. This yields the Markov chain $\{(D_n, N_n)\}$ on $[0, \infty) \times \{1, 2, \ldots\}$ with transitions controlled by

$$\mathbb{P}[D_{n+1} = D_n, N_{n+1} = N_n - 1 | D_n, N_n] = 1 \quad \text{if } N_n \geq 2;$$

$$\mathbb{P}[D_{n+1} \in E | D_n = z, N_n = 1] = \mathbb{P}[Y_1 \in E | Y_0 = z]$$

for all measurable $E \subseteq [1, \infty)$;

$$\mathbb{P}[N_{n+1} = S(D_{n+1}) | D_n, N_n = 1, D_{n+1}] = 1.$$

Using the standard equilibrium equations, if $\tilde{\pi}$ is the equilibrium distribution of $(D, N)$, then

$$\tilde{\pi}(z, 1) = \tilde{\pi}(z, 2) = \cdots = \tilde{\pi}(z, S(z))$$

and thus $\pi_D(z) = \tilde{\pi}(z, \cdot) \propto \pi_Y(z)S(z)$. Hence, the equilibrium distribution of $D$ is the equilibrium of $Y$ reweighted using $S$. It is a classical probability result [10] that under stationarity the number of people in the $D/M/1$ queue (used in the construction of $Y$) is geometric with parameter $\eta$, where $\eta$ is the smallest positive root of

$$\eta = \beta^{1-\eta}.$$

(Note that $0 < \eta < 1$ since $\beta < e^{-1}$.) Thus the equilibrium distribution of the queue workload $U$ is exponential of rate $(1 - \eta)$. Since $Y \propto \exp(U)$, the equilibrium density of $Y$, $\pi_Y$, satisfies

(14) $$\pi_Y(z) \propto z^{-(2-\eta)}.$$

Reweighting $Y$ using $S$ yields the equilibrium density of $D$,

(15) $$\pi_D(z) \propto S(z) z^{-(2-\eta)}.$$

A suitable pause function $S$ must therefore satisfy $S(z) < z^{1-\eta}$ in order to obtain a probability density in (15). The dominating process constructed in the proof of Theorem 16 requires $F \leq S$ and hence this imposes the restriction $F(z) < z^{1-\eta}$; in particular, this means that $F(z)/z \to 0$ as $z \to \infty$.

3.2. *Tame and wild chains.* The above discussion motivates the following definition of a *tame* chain. We write $\lceil z \rceil$ to denote the smallest integer greater than or equal to $z$.

DEFINITION 14. A Markov chain $X$ is *tame with respect to a scale function* $V$ if the following two conditions hold:

(a) there exists a small set $C' := \{x : V(x) \leq d'\}$ and a nondecreasing *taming function* $F : [1, \infty) \to \{1, 2, \ldots\}$ of the form

(16) $$F(z) = \begin{cases} \lceil \lambda z^\delta \rceil, & z > d', \\ 1, & z \leq d', \end{cases}$$

for some constants $\lambda > 0$, $\delta \in [0, 1)$ such that the chain $X' = \{X_{\theta_n}\}$ possesses the drift condition

(17) $$PV \leq \beta V + b' \mathbf{1}_{C'},$$

where $\{\theta_n\}$ is an adaptive sampling scheme defined using $F$, as in (12);

(b) the constant $\beta$ in inequality (17) satisfies

(18) $$\log \beta < \delta^{-1} \log(1 - \delta).$$

We say that $X$ is *tamed* (*with respect to* $V$) *by the function* $F$. We may simply say that $X$ is *tame*, without mention of a specific scale function. A chain that is not tame is said to be *wild*.

Thus a tame chain is one for which we can exhibit an explicit adaptive subsampling scheme using a power function $F$ and for which the subsampled chain so produced is geometrically ergodic with sufficiently small $\beta$.

Note that all geometrically ergodic chains are trivially tame: if $X$ satisfies Condition $\mathrm{GE}(V, \beta, b, C)$, then $X$ is tamed by the function

$$F(z) = k \qquad \text{for } z > \sup_{y \in C} V(y),$$

for any integer $k > 1 + 1/\log \beta$.

Definition 14 is strongly motivated by the discussion in Section 3.1. From (16), we see that $F$ produces a simple adaptive subsampling scheme, as in Definition 13. $F$ is also a nondecreasing function, which accords with our intuition; if $V(X)$ is large, then we expect to wait longer before subsampling again, to create enough drift in the chain to produce a geometric Foster–Lyapunov condition. Requirement (b) of Definition 14 is made for two reasons. First, it ensures that $\beta < e^{-1}$ and so ensures ergodicity of the $D/M/1$ queue workload $U$ used in the construction of $Y$. Second, it ensures that the weighted equilibrium distribution of $Y$ using $S$ (as described at the end of Section 3.1) is a proper distribution; this will be shown in the proof of Theorem 16.

Kendall [14] shows that a dominating process exists for $V(X')$ even if $\beta > e^{-1}$, but recall that this involves a further subsampling of $X'$ with a fixed period $k$. Here, $\beta < e^{-1}$ is made a requirement of the adaptive subsampling process to avoid this situation, since further subsampling of $X'$ would result in a composite nondeterministic subsampling scheme.

The main theorem of this paper is the following.

THEOREM 15. *Suppose that $X$ is tame with respect to a scale function $V$. Then there exists a domCFTP algorithm for $X$ using a dominating process based on $V$.*

Theorem 15 is true for all geometrically ergodic chains by the result of [14]. As with the results of [7] and [14], this algorithm may not be implementable in practice. The proof of Theorem 15 results directly from Theorem 16 and the discussion in Section 3.3 below, where a description of the *domCFTP* algorithm is given.

THEOREM 16. *Suppose that $X$ satisfies the weak drift condition $PV \leq V + b\mathbf{1}_C$ and that $X$ is tamed with respect to $V$ by the function*

$$F(z) = \begin{cases} \lceil \lambda z^\delta \rceil, & z > d', \\ 1, & z \leq d', \end{cases}$$

*with the resulting subsampled chain $X'$ satisfying a drift condition $PV \leq \beta V + b' \mathbf{1}_{[V \leq d']}$, with $\log \beta < \delta^{-1} \log(1-\delta)$. Then there exists a stationary ergodic process $D$ which dominates $V(X)$ at the times $\{\sigma_n\}$ when $D$ moves.*

PROOF. We shall construct a Markov chain $(D, N)$ by starting with a process $Y$ and pausing it using a function $S$, to be defined shortly. Before beginning the main calculation of the proof, we define some constant. These are determined explicitly from the taming function $F$ and the drift conditions satisfied by $X$ and $X'$. First, choose $\beta^* > \beta$ such that

(19) $$\log \beta < \log \beta^* < \delta^{-1} \log(1-\delta).$$

(That this is possible is a result of the definition of tameness.) Then set

$$\begin{aligned} a &= \frac{b'}{1-\beta}(1 + b(\lambda+1)), \\ d^* &= \min\{z \geq d' : (\beta^* - \beta)z \geq b(\lambda+1)z^\delta + a\}, \\ b^* &= b(\lambda+1)d^{*\delta} + a, \\ h^* &= d^* + \frac{b^*}{\beta^*}. \end{aligned}$$

Finally, consider the set $C^* = \{x : V(x) \leq h^*\}$. As a sublevel set, $C^*$ is $m$-small for some integer $m \geq 1$. We are now in a position to define the function $S$:

$$S(z) = \begin{cases} (m \vee F(h^*)) \lceil \lambda z^\delta \rceil, & z \geq h^*, \\ (m \vee F(h^*)), & z < h^*. \end{cases} \tag{20}$$

Note that $F(x) \leq S(z)$ for all $x \leq z$ (since $h^* \geq d'$) and that

$$S(z) \geq m \vee F(h^*) \geq m \quad \text{for all } z \geq 0. \tag{21}$$

Define the process $Y = h^* \exp(U)$, where $U$ is the system workload of a $D/M/1$ queue with arrivals every $\log(1/\beta^*)$ time units and service times being independent and of unit exponential distribution. Positive recurrence of $U$ follows from (19). Pause $Y$ using $S$ (as described on page 9) and call the resulting process $D$. The stationary distribution of $D$, as shown at the end of Section 3.1, is given by

$$\begin{aligned} \pi_D(z) &\propto S(z) z^{-(2-\eta)} \\ &\propto z^{-(2-\eta-\delta)} \quad \text{(for } z > h^*\text{)}, \end{aligned} \tag{22}$$

where $\eta < 1$ is the smallest positive solution to the equation

$$\eta = \beta^{*(1-\eta)}.$$

Now, by our choice of $\beta^*$ above, we have

$$(1-\eta)^{-1} \log \eta = \log \beta^* < \delta^{-1} \log(1-\delta),$$

so $\eta < 1 - \delta$. Hence, $2 - \eta - \delta > 1$, so we see from (22) that $\pi_D$ is a proper density.

Suppose that $(D_{\sigma_n}, N_{\sigma_n}) = (z, S(z))$ and that $V(X_{\sigma_n}) = V(x) \leq z$. We wish to show that $D_{\sigma_{n+1}}$ dominates $V(X_{\sigma_{n+1}})$, where $\sigma_{n+1} = \sigma_n + S(z)$ is the time at which $D$ next moves. Domination at successive times $\{\sigma_j\}$ at which $D$ moves then follows inductively. For simplicity in the calculations below, we set $\sigma_n = 0$.

Let $\{\theta_n\}$ be the adaptive subsampling scheme for $X$ defined recursively by the taming function $F$. Define a region $R(z) \subset \mathcal{X} \times \mathbb{Z}_+$ to be the so-called "short sampling" region:

$$R(z) = \{(y, t) : F(V(y)) + t > S(z)\}.$$

In other words, once the chain $\{X_{\theta_n}, \theta_n\}$ hits the (deterministic) region $R(z)$ (at time $\theta_j$, say), the next subsampling time $[\theta_{j+1} = \theta_j + F(V(X_{\theta_j}))]$ will lie beyond the time $S(z)$ at which the dominating process moves (see Figure 1). Define

$$T(z) = \min\{\theta_n : (X_{\theta_n}, \theta_n) \in R(z)\}$$

to be a stopping time for $X$ and define

$$T'(z) = \min\{n : (X_{\theta_n}, \theta_n) \in R(z)\}$$

to be the associated stopping time for $X'$. [Note that $T'(z) \geq 1$ since $V(x) \leq$

$z$ implies that $F(V(x)) \leq S(z)$.]

Our aim is to control $\mathbb{E}_x[V(X_{S(z)})]$, recalling that $V \geq 1$ and that $S(z)$ is deterministic:

$$\mathbb{E}_x[V(X_{S(z)})] = \mathbb{E}_x[\mathbb{E}_{X_{T(z)}}[V(X_{S(z)})]]$$

$$\leq \mathbb{E}_x\left[V(X_{T(z)}) + b\sum_{j=T(z)}^{S(z)-1} \mathbb{P}_{X_{T(z)}}[X_j \in C]\right]$$

using the weak drift condition of the theorem

$$\leq \mathbb{E}_x[V(X_{T(z)})] + b\mathbb{E}_x[(S(z) - T(z))]$$

(23)
$$\leq \mathbb{E}_x[V(X_{T(z)})] + b\mathbb{E}_x[F(V(X_{T(z)}))]$$

since $S(z) - T(z) < F(V(X_{T(z)}))$,

by the definition of $R(z)$

$$\leq \mathbb{E}_x[V(X_{T(z)})] + b(\lambda + 1)\mathbb{E}_x[V(X_{T(z)})^\delta]$$

by the definition of $F$.

Now, the chain $X' = \{X_{\theta_n}\}$ is geometrically ergodic (since $X$ is tamed by $F$), so Lemma 6 tells us that

(24)
$$\mathbb{E}_x[V(X_{T(z)})] = \mathbb{E}_x[V(X'_{T'(z)})] \leq \beta V(x) + \frac{b'}{1-\beta}.$$

Furthermore, Lemma 5 yields

$$\mathbb{E}_x[V(X_{T(z)})^\delta] = \mathbb{E}_x[V(X'_{T'(z)})^\delta]$$

(25)
$$\leq \beta^\delta V^\delta(x) + \left(\frac{b'}{1-\beta}\right)^\delta \mathbf{1}_{[V(x)\leq d']}$$

$$\leq V^\delta(x) + \frac{b'}{1-\beta}.$$

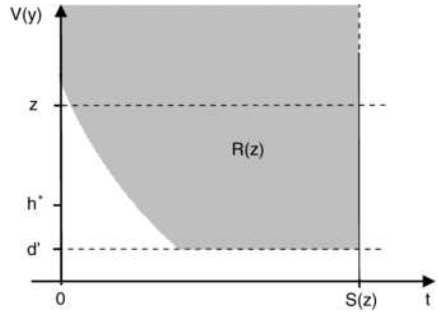

FIG. 1. *Depiction of the region $R(z)$.*

Combining equations (23), (24) and (25) and making use of the constants defined at the start of this proof, we obtain

$$\mathbb{E}_x[V(X_{S(z)})] \leq \beta V(x) + b(\lambda+1)V^\delta(x) + a$$
(26)
$$\leq \beta^* V(x) + b^* \mathbf{1}_{[V(x) \leq d^*]}.$$

Thus a geometric drift condition holds at time $S(z)$ for all chains $V(X)$ with starting states $x$ satisfying $V(x) \leq z$. As in the proof of Theorem 12, it follows from inequality (26) that $V(X_{S(z)})$ can be dominated by $D_{S(z)}$ [17]. □

Note that questions 1 and 3 at the start of Section 3 have now been answered: we have defined what is meant by an adaptive subsampling scheme and shown that if this takes a particular (power function) form then a stationary process $D$ that dominates $V(X)$ at times $\{\sigma_n\}$ can be produced.

3.3. *The domCFTP algorithm for tame chains.* In this section, we describe the *domCFTP* algorithm for tame chains and hence complete the proof of Theorem 15. We begin this by answering question 4 of page 8, by showing how to simulate $(D, N)$ in equilibrium and in reversed time. Furthermore, this simulation is quite simple to implement when the function $S$ is of the form (20).

The first point to make here is that one can easily simulate from $\pi_D$ using rejection sampling [21]: using (15), for some constant $\gamma > 0$, we have

$$\pi_D(z) = \gamma \left( \frac{1}{2} \frac{\lceil \lambda z^\delta \rceil}{\lambda z^\delta} \right) \frac{1}{z^{2-\eta-\delta}}$$
$$= \gamma p(z) g(z),$$

where $p(z) \in [1/2, 1]$ and $g(z)$ is a Pareto density (since $2 - \eta - \delta > 1$, as in the proof of Theorem 16). Now, given $D_0 = z_0$ as a draw from $\pi_D$, set $N_0 := n_0$, where $n_0 \sim \text{Uniform}\{1, 2, \ldots, S(z_0)\}$. It follows from the construction of $(D, N)$ in Section 3.1 that $(D_0, N_0) \sim \tilde{\pi}$, as required.

The chain $(D, N)$ is then simple to run in reversed time using the facts that the jumps of $D$ are those of the underlying exponential queue workload process $Y$ and that the pause function $S$ is deterministic. (Recall the forward construction on page 9 and see Figure 2. More details can be found in [3].)

We now show that $D$ is a dominating process for $X$ (at the times when $D$ moves) based on the scale function $V$, with threshold $h^*$ (recall Definition 11). Also, recall from the proof of Theorem 16 that the set $C^* = \{x : V(x) \leq h^*\}$ is $m$-small.

First, the proof of Theorem 16 shows that the link between stochastic domination and coupling [17] may be exploited to couple the various $X^{x,\sigma_{-M}}$ with $D$ such that for all $n \leq M$,

(27) $\quad V(X^{x,\sigma_{-M}}_{\sigma_{-n}}) \leq D_{\sigma_{-n}} \implies V(X^{x,\sigma_{-M}}_{\sigma_{-(n-1)}}) \leq D_{\sigma_{-(n-1)}}.$

We now turn to part (b) of Definition 11. Since $C^*$ is $m$-small, there exists a probability measure $\nu$ and a scalar $\varepsilon \in (0,1)$ such that for all Borel sets $B \subset [1,\infty)$, whenever $V(x) \leq h^*$,

$$\mathbb{P}[V(X_m) \in B | X_0 = x] \geq \varepsilon \nu(B).$$

Therefore, since $S(h^*) \geq m$ [as noted in (21)],

$$\mathbb{P}[V(X_{S(h^*)}) \in B | X_0 = x] \geq \varepsilon P_\nu^{S(h^*)-m}(B),$$

so $C^*$ is $S(h^*)$-small. Furthermore, the stochastic domination which has been arranged in the construction of $D$ means that for all $u \geq 1$, whenever $V(x) \leq y$,

$$[V(X_{S(y)}) > u | X_0 = x] \leq \mathbb{P}[Y_1 > u | Y_0 = y].$$

We can couple in order to arrange for regeneration if a probability measure $\tilde{\nu}$ can be identified, defined solely in terms of $P_\nu^{S(h^*)-m}$ and the dominating jump distribution $\mathbb{P}[Y_1 \geq u | Y_0 = y]$, such that for all $u \geq 1$, whenever $V(x) \leq y$,

$$\mathbb{P}[V(X_{S(y)}) > u | X_0 = x] - \varepsilon P_\nu^{S(h^*)-m}((u,\infty))$$
$$\leq \mathbb{P}[Y_1 > u | Y_0 = y] - \varepsilon \tilde{\nu}((u,\infty))$$

$$P_\nu^{S(h^*)-m}((u,\infty)) \leq \tilde{\nu}((u,\infty));$$

and

$$\mathbb{P}[Y_1 \in E | Y_0 = y] \geq \varepsilon \tilde{\nu}(E)$$

for all measurable $E \subseteq [1,\infty)$.

Recall the following result, a proof of which is provided in [14].

LEMMA 17. *Suppose that $U$, $V$ are two random variables defined on $[1,\infty)$ such that*:

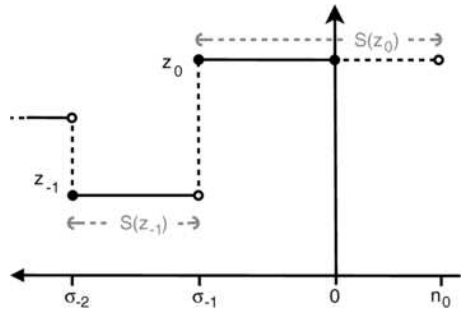

FIG. 2. *Construction of $D$ in reversed time.*

(a) *The distribution $\mathcal{L}(U)$ is stochastically dominated by the distribution $\mathcal{L}(V)$, that is,*

$$\mathbb{P}[U > u] \leq \mathbb{P}[V > u] \qquad \text{for all } u \geq 1;$$

(b) *$U$ satisfies a minorization condition: for some $\beta \in (0,1)$ and probability measure $\psi$,*

$$\mathbb{P}[U \in E] \geq \beta \psi(E) \qquad \text{for all Borel sets } E \subseteq [1, \infty).$$

*Then there exists a probability measure $\mu$ stochastically dominating $\psi$ and such that $\beta\mu$ is minorized by $\mathcal{L}(V)$. Moreover, $\mu$ depends only on $\beta\psi$ and $\mathcal{L}(V)$.*

Therefore, using Lemma 17, $\mathcal{L}(X_{\sigma_{-(n-1)}} | X_{\sigma_{-n}} = x)$ may be coupled to $\mathcal{L}(D_{\sigma_{-(n-1)}} | D_{\sigma_{-n}} = y)$ whenever $V(x) \leq y$, in a way that implements stochastic domination and ensures that all of the $X_{\sigma_{-(n-1)}}$ can regenerate simultaneously whenever $D_{\sigma_{-n}} \leq h^*$.

Finally, it is easy to see that part (c) of Definition 11 is satisfied: the system workload $U$ of the queue will hit zero infinitely often and therefore $D$ will hit level $h^*$ infinitely often.

We can now describe a *domCFTP* algorithm based on $X$ which yields a draw from the equilibrium distribution.

ALGORITHM.

- Simulate $D$, as a component of the stationary process $(D, N)$, backward in time until the most recent $\sigma_{-M} < 0$ for which $D_{\sigma_{-M}} \leq h^*$;
- while coalescence does not occur at time $\sigma_{-M}$, extend $D$ backward until the most recent $\sigma_{-M'} < \sigma_{-M}$ for which $D_{\sigma_{-M'}} \leq h^*$ and set $M \leftarrow M'$;
- starting with the unique state produced by the coalescence event at time $\sigma_{-M}$ simulate the coupled $X$ forward at times $\sigma_{-M}, \sigma_{-(M-1)}, \sigma_{-(M-2)}, \ldots$, up to and including time $\sigma_{-1}$;
- run the chain $X$ forward (from its unique state) from time $\sigma_{-1}$ to $0$ (see Figure 3);
- return $X_0$ as a perfect draw from equilibrium.

LEMMA 18. *The output of the above algorithm is a draw from the stationary distribution of the target chain $X$.*

PROOF. The stochastic domination of (27) and Theorem 2.4 of [17], Chapter IV guarantee the existence of a joint transition kernel $P_{X,D}^{(n)}$ that provides domination of $X$ by $D$ and such that the marginal distributions of

$X$ and $D$ are correct. That is, for $x \leq y$, with $n = S(y)$, for all $z \geq 1$,

$$P^{(n)}_{X,D}(x,y; V^{-1}((z,\infty)), [1,z]) = 0,$$

$$\int_{V^{-1}([1,z])} \int_1^\infty P^{(n)}_{X,D}(x,y; du, dv) = P^{(n)}_X(x; V^{-1}([1,z])),$$

$$\int_\mathcal{X} \int_1^z P^{(n)}_{X,D}(x,y; du, dv) = P^{(n)}_D(y; [1,z]).$$

The chains $X$ and $D$ (run forward) may therefore be constructed in either of the following two ways.

1. Given $D_{\sigma_{-m}}$ and $X_{\sigma_{-m}} \leq D_{\sigma_{-m}}$, with $n = S(D_{\sigma_{-m}})$:
   - draw $D_{\sigma_{-(m-1)}}$ from the probability kernel
   $$P^{(n)}_D(D_{\sigma_{-m}}; \cdot);$$
   - draw $X_{\sigma_{-(m-1)}}$ from the regular conditional probability
   $$\frac{P^{(n)}_{X,D}(X_{\sigma_{-m}}, D_{\sigma_{-m}}; \cdot, D_{\sigma_{-(m-1)}})}{P^{(n)}_D(D_{\sigma_{-m}}; D_{\sigma_{-(m-1)}})};$$
   - draw $X_{\sigma_{-m}+1}, X_{\sigma_{-m}+2}, \ldots, X_{\sigma_{-(m-1)}-1}$ as a realization of $X$ conditioned on the values of $X_{\sigma_{-m}}$ and $X_{\sigma_{-(m-1)}}$ (i.e., as a *Markov bridge* between $X_{\sigma_{-m}}$ and $X_{\sigma_{-(m-1)}}$).

2. Given $D_{\sigma_{-m}}$ and $X_{\sigma_{-m}} \leq D_{\sigma_{-m}}$, with $n = S(D_{\sigma_{-m}})$:
   - draw $X_{\sigma_{-m}+1}, X_{\sigma_{-m}+2}, \ldots, X_{\sigma_{-(m-1)}}$ using the normal transition kernel for $X$, noting that the distribution of $X_{\sigma_{-(m-1)}}$ is exactly the same as if it were drawn directly from $P^{(n)}_X(X_{\sigma_{-m}}; \cdot)$;

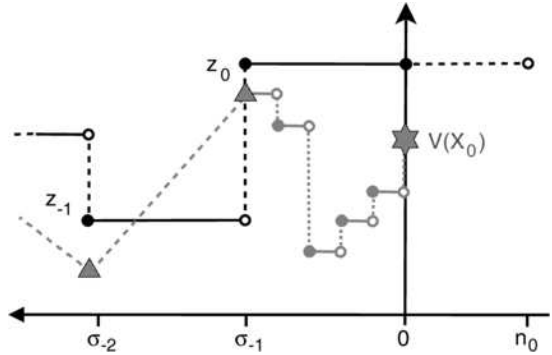

FIG. 3. *Final stage of the domCFTP algorithm: $D$ (black circles •) dominates $V(X)$ (red triangles ▲) at times $\{\sigma_n\}$. To obtain the draw from equilibrium, $X_0$, $X$ can be run from time $\sigma_{-1}$ to 0 without reference to $D$ after time $\sigma_{-1}$.*

- draw $D_{\sigma_{-(m-1)}}$ from the regular conditional probability

$$P^{(n)}_{D|\{X\}}(\cdot|D_{\sigma_{-m}}, X_{\sigma_{-m}}, X_{\sigma_{-m}+1}, \ldots, X_{\sigma_{-(m-1)}})$$
$$= \frac{P^{(n)}_{X,D}(X_{\sigma_{-m}}, D_{\sigma_{-m}}; X_{\sigma_{-(m-1)}}, \cdot)}{P^{(n)}_X(X_{\sigma_{-m}}; X_{\sigma_{-(m-1)}})}.$$

Each of these two methods produces chains $X$ and $D$ which satisfy the stochastic domination of (27). Method 1 is that which is effectively used by the algorithm, although there is no need for the final superfluous step (the Markov bridge) when implementing the algorithm. Method 2, however, makes it clear that $X$ has the correct Markov transition kernel to be the required target chain. Furthermore, the equivalence of the two schemes proves the validity of the final step of the algorithm, where the chain $X$ is run from time $\sigma_{-1}$ to 0 without reference to $D$.

Finally, the proof that the algorithm returns a draw from equilibrium follows a standard renewal theory argument. Consider a stationary version of the chain $X$, say $\hat{X}$, run from time $-\infty$ to 0. The regenerations of $\hat{X}$ (when it visits the small set $C^*$) and those of $D$ (when it hits level $h^*$) form two positive recurrent renewal processes (with that of $\hat{X}$ being aperiodic). Therefore, if $D$ is started far enough in the past, then there will be a time $-T$ at which both $\hat{X}$ and $D$ regenerate simultaneously. Now, consider the process $\tilde{X}_n = \hat{X}_n \mathbf{1}_{[n<-T]} + X_n \mathbf{1}_{[n\geq -T]}$. Clearly, $\tilde{X}$ is stationary and follows the same transitions of $X$ from time $-T$ to 0. Thus $X_0 = \tilde{X}_0 \sim \pi$, so the output as the algorithm is indeed a draw from the required equilibrium distribution. □

This concludes the proof of Theorem 15. We have produced a *domCFTP* algorithm based on the scale function $V$ for the tame chain $X$.

3.4. *When is a chain tame?* As a consequence of Theorem 15, question 2 of page 8 can be rephrased as: when is a chain tame? Note that a tame chain will not necessarily be tamable with respect to all scale functions, of course.

In this section, we present an equivalent definition of tameness and prove some sufficient conditions for a polynomially ergodic chain to be tame. The following theorem shows that tameness is determined precisely by the behavior of the chain until the time that it first hits the small set $C$.

THEOREM 19. *Suppose that $X$ satisfies the weak drift condition $PV \leq V + b\mathbf{1}_C$. Then for $n(x) = o(V(x))$, the following two conditions are equivalent:*

(i) *there exists $\beta \in (0,1)$ such that $\mathbb{E}_x[V(X_{n(x)})] \leq \beta V(x)$ for $V(x)$ sufficiently large;*

(ii) *there exists $\beta' \in (0,1)$ such that $\mathbb{E}_x[V(X_{n(x) \wedge \tau_C})] \leq \beta' V(x)$ for $V(x)$ sufficiently large.*

*Furthermore, if $V(x)$ is sufficiently large, we may take $|\beta - \beta'| < \varepsilon$ for any $\varepsilon > 0$.*

PROOF. Since $C = \{x : V(x) \leq d\}$ is a sublevel set, we can split the expectation of $V(X_{n(x) \wedge \tau_C})$ according to whether $\tau_C \leq n(x)$ or not, to show that

$$\mathbb{E}_x[V(X_{n(x) \wedge \tau_C})] \leq \sup_{y \in C} V(y) + \mathbb{E}_x[V(X_{n(x)}); \tau_C > n(x)]$$

$$\leq \sup_{y \in C} V(y) + \mathbb{E}_x[V(X_{n(x)})],$$

so (i) $\Rightarrow$ (ii).

We now prove the reverse implication. Using the weak drift condition for $X$ and recalling that $n(x)$ is deterministic, we have

$$\mathbb{E}_x[V(X_{n(x)}); \tau_C \leq n(x)] = \sum_{k=1}^{n(x)} \mathbb{E}_x[\mathbb{E}_{X_k}[V(X_{n(x)-k})]; \tau_C = k]$$

$$\leq \sum_{k=1}^{n(x)} \sup_{y \in C} \mathbb{E}_y[V(X_{n(x)-k}) | X_k = y] \mathbb{P}_x[\tau_C = k]$$

$$\leq \sum_{k=1}^{n(x)} \sup_{y \in C}(V(y) + b(n(x) - k)) \mathbb{P}_x[\tau_C = k]$$

$$\leq d + n(x)b.$$

Assuming (ii), we therefore have

$$\mathbb{E}_x[V(X_{n(x)})] \leq \mathbb{E}_x[V(X_{n(x) \wedge \tau_C})] + \mathbb{E}_x[V(X_{n(x)}); \tau_C \leq n(x)]$$

$$\leq \beta' V(x) + d + n(x)b$$

$$\leq \beta V(x)$$

for all sufficiently large $V(x)$, since $n(x) = o(V(x))$.

Finally, due to the restriction on the size of $n(x)$, it is clear that $\beta$ and $\beta'$ may be made arbitrarily close by simply restricting attention to $x$ for sufficiently large $V(x)$. □

Suppose that we now modify the behavior of a tame chain $X$ when it is in the small set $C$. The following simple corollary of Theorem 19 shows that provided the resulting chain still satisfies a weak drift condition, tameness is preserved under such modification.

COROLLARY 20. *Suppose $X$ satisfies the drift condition $PV \leq V + b\mathbf{1}_C$ and that $X$ is tamed by the function $F$ to produce a chain $X'$ satisfying $GE(V, \beta, b', C')$. Let $\hat{X}$ be a new chain produced by modifying the behavior of $X$ when in $C$, such that $\hat{X}$ satisfies $PV \leq V + \hat{b}\mathbf{1}_C$. Then $F$ also tames $\hat{X}$, and the resulting chain $\hat{X}'$ satisfies $GE(V, \hat{\beta}, \hat{b}', \hat{C}')$ for any $\hat{\beta}' \in (\beta, 1)$.*

PROOF. Write $F_x = F(V(x))$. Since $X$ is tame, Theorem 19 tells us that for $V(x)$ sufficiently large,
$$\mathbb{E}_x[V(X_{F_x \wedge \tau_C})] \leq \tilde{\beta} V(x)$$
for any $\tilde{\beta} \in (\beta, 1)$. Now, since
$$\hat{X} \mathbf{1}_{[\hat{\tau}_C \geq F_x]} \stackrel{d}{=} X \mathbf{1}_{[\tau_C \geq F_x]},$$
by definition,
$$\mathbb{E}_x[V(\hat{X}_{F_x \wedge \hat{\tau}_C})] \leq \tilde{\beta} V(x).$$
Furthermore, since $\hat{X}$ satisfies the drift condition $PV \leq V + \hat{b}\mathbf{1}_C$, a second application of Theorem 19 yields
$$\mathbb{E}_x[V(\hat{X}_{F_x})] \leq \hat{\beta} V(x),$$
where $\hat{\beta} \in (\tilde{\beta}, 1)$ may be chosen arbitrarily close to $\tilde{\beta}$ (and hence to $\beta$). Thus the same function $F$ also tames $\hat{X}$. □

We have already remarked that all geometrically ergodic chains are tame. The next two theorems provide sufficient conditions for a polynomially ergodic chain to be tame.

THEOREM 21. *Let $X$ be a chain satisfying a drift condition $PV \leq V - cV^\alpha + b\mathbf{1}_C$ for which $V(X)$ has bounded upward jumps whenever $X \notin C$. That is, $V(X_1) \leq V(X_0) + K$ whenever $X_0 \notin C$, for some constant $K < \infty$. Then $X$ is tame.*

PROOF. From Theorem 19, we see that it is sufficient to show that by choosing an appropriate taming function $F$, we can obtain the bound

(28) $\quad \mathbb{E}_x[V(X_{F(V(x))}); F(V(x)) < \tau_C] \leq \beta V(x) + b' \mathbf{1}_{C'}(x).$

Choose $\beta$ sufficiently small to satisfy

(29) $\quad\quad\quad\quad\quad \log \beta < (1-\alpha)^{-1} \log \alpha$

and then choose $\lambda$ sufficiently large so that $\lambda^{-1} < \beta c(1-\alpha)$. Define the constant $d_1$ by

$$d_1^\alpha = \frac{K}{c(1-\alpha)\lambda}\left(\beta - \frac{1}{c(1-\alpha)\lambda}\right)^{-1}$$

and define $C_1 = \{x : V(x) \leq d_1\}$. Note that if $x \notin C_1$, then

(30) $$\left(\beta - \frac{1}{c(1-\alpha)\lambda}\right)V(x) \geq \left(\frac{K}{c(1-\alpha)\lambda}\right)V^{1-\alpha}(x).$$

Finally, set $d' = \max\{d, d_1\}$ and let $C' = \{x : V(x) \leq d'\}$.

Now, define the taming function $F$ by

(31) $$F(z) = \begin{cases} \lceil \lambda z^{1-\alpha} \rceil, & \text{for } z > d', \\ 1, & \text{for } z \leq d'. \end{cases}$$

Write $F_x = F(V(x))$ to ease notation. Then for $x \notin C'$, since the upward jumps of $V(X)$ before time $\tau_C$ are bounded above by $K$, we have

$$\mathbb{E}_x[V(X_{F_x}); F_x < \tau_C] \leq (V(x) + KF_x)\mathbb{P}_x[\tau_C > F_x]$$

$$\leq (V(x) + KF_x)\frac{\mathbb{E}_x[\tau_C]}{F_x} \quad \text{by Markov's inequality}$$

$$\leq (V(x) + KF_x)\frac{V^{1-\alpha}(x)}{c(1-\alpha)F_x} \quad \text{by Corollary 9}$$

$$\leq \frac{V(x)}{c(1-\alpha)\lambda} + \left(\frac{K}{c(1-\alpha)\lambda}\right)V^{1-\alpha}(x) \quad \text{using (31)}$$

$$\leq \beta V(x) \quad \text{by inequality (30)}.$$

Finally, for $x \in C'$, we have

$$\mathbb{E}_x[V(X_{F_x})] = \mathbb{E}_x[V(X_1)] \leq V(x) + b$$

$$\leq \beta V(x) + (1-\beta)d' + b$$

$$= \beta V(x) + b',$$

where $b' = (1-\beta)d' + b < \infty$. Hence, (28) is satisfied for all $x$ and $X$ is tame. □

The following proof makes use of Proposition 10, which was borrowed from [11]. Note that tameness is clearly monotonic in the drift exponent $\alpha$ since chains satisfying $PE(V, c, \alpha, b, C)$ also satisfy $PE(V, c, \alpha', b, C)$ for all $\alpha' \leq \alpha$.

THEOREM 22. *Let $X$ be a chain satisfying the drift condition $PV \leq V - cV^\alpha + b\mathbf{1}_C$, with $\alpha > 3/4$. Then $X$ is tame.*

PROOF. Let $\rho = (1-\alpha)^{-1}/2 > 2$ and set $\alpha' = 2\alpha - 1$. Writing $V_\rho = V^{1-\rho(1-\alpha)} = V^{1/2}$ and using Lemma 8, we have

$$PV_\rho \leq V_\rho - V_\rho^{\alpha'} + b_1 \mathbf{1}_C$$

for some $b_1 < \infty$. We shall seek a time change that produces a geometric Foster–Lyapunov condition on this scale, $V_\rho$. As in the proof of Theorem 21, we simply need to control

$$\mathbb{E}_x[V_\rho(X_{F_x}); F_x < \tau_C],$$

where $F_x = F(V_\rho(x))$.

By Proposition 10,

$$\mathbb{C}E_x\left[\sum_{n=0}^{\tau_C - 1} n^{\rho-1} V_\rho(X_n)\right] \leq MV(x)$$

for some constant $M < \infty$. Thus

(32) $$\mathbb{E}_x[V_\rho(X_{F_x}); F_x < \tau_C] \leq \frac{MV(x)}{F_x^{\rho-1}}.$$

Now, choose $\beta > 0$ such that $\log \beta < (\rho - 1)\log((\rho-2)/(\rho-1))$ and define the taming function $F$ by

$$F(z) = \lceil (\lambda z)^{1/(\rho-1)} \rceil \vee 1$$

for any $\lambda > M/\beta$. Then from inequality (32),

$$\mathbb{E}_x[V_\rho(X_{F_x})] \leq \frac{MV(x)}{F_x^{\rho-1}}$$
$$\leq \beta V_\rho(x)$$

for $V_\rho(x)$ sufficiently large. Therefore, $F$ tames $X$, as required. $\square$

In fact, it turns out that *any* chain satisfying drift Condition PE may be adaptively subsampled as above to produce a geometrically ergodic chain (see [2] for details). However, for $\alpha \leq 3/4$, the pause function produced leads to an improper equilibrium distribution for the dominating process of Theorem 16. Connor [2] shows how this lower bound on $\alpha$ may be further reduced to 0.704, but tameness for $\alpha \leq 0.704$ remains to be proven. This is not to say, of course, whether or not there may exist another suitable pause function, possibly on a different scale.

These two sufficient conditions are not necessary for a chain to be tame: in Section 4.4, we present an example of a chain that satisfies Condition PE with drift coefficient $\alpha = 1/2$ and which does not have bounded jumps for $X \notin C$, and we show explicitly that it is tame.

**4. Examples.** We now present four explicit examples of polynomially ergodic chains and show that they are tame. The first two of these are tame by Theorem 21 and the third by Theorem 22. The final example, (4.4), shows that the sufficient conditions of Theorems 21 and 22 are not necessary for $X$ to be tame.

4.1. *Epoch chain.* Consider the Markov chain $X$ on $\{0, 1, 2, \ldots\}$ with the following transition kernel: for all $x \in \{0, 1, 2, \ldots\}$,

$$P(x, x) = \theta_x; \qquad P(0, x) = \zeta_x;$$
$$P(x, 0) = 1 - \theta_x.$$

Thus $X$ spends a random length of time (an epoch) at level $x$ before jumping to 0 and regenerating. Meyn and Tweedie ([18], page 362) show that this chain is ergodic if $\zeta_x > 0$ for all $x$ and

(33) $$\sum_x \zeta_x (1 - \theta_x)^{-1} < \infty.$$

Furthermore, they show that the chain is not geometrically ergodic if $\theta_x \to 1$ as $x \to \infty$, regardless of how fast $\zeta_x \to 0$.

Now, suppose that $\theta_x = 1 - \kappa(x+1)^{-\lambda}$ for some suitable $\kappa, \lambda > 0$. We now slightly strengthen condition (33) on $\{\zeta_x\}$ to obtain a polynomial drift condition: we require that there exists $\varepsilon > 0$ such that $\sum_x \zeta_x x^{(1+\varepsilon)\lambda} < \infty$.

Let $C = [0, \kappa^{1/\lambda}]$. Then following drift condition holds:

(34) $$\mathbb{E}_x[V(X_1)] \leq V(x) - \kappa V^\alpha(x) + b\mathbf{1}_C(x),$$

where $V(x) = (x+1)^m$, $m = (1+\varepsilon)\lambda$ and $\alpha = \varepsilon/(1+\varepsilon)$. This chain then satisfies the conditions of Theorem 21 and is therefore tame.

4.2. *Delayed death process.* Consider the Markov chain $X$ on $\{0, 1, 2, \ldots\}$ with the following transition kernel:

$$P(x, x) = \theta_x, \qquad x \geq 1,$$
$$P(x, x-1) = 1 - \theta_x, \qquad x \geq 1,$$
$$P(0, x) = \zeta_x > 0, \qquad x \in \{0, 1, 2, \ldots\},$$

where $\theta_x = 1 - \kappa(x+1)^{-\lambda}$ for some suitable $\kappa > 0, \lambda > 1$, and $\zeta_x \to 0$ as $x \to \infty$ sufficiently fast to ensure that

$$\mathbb{E}_0[\tau_0] = 1 + \sum_{x=1}^{\infty} \zeta_x \sum_{y=1}^{x} (1 - \theta_y)^{-1} < \infty,$$

making $X$ ergodic.

It is simple to show that $X$ is not geometrically ergodic, but that it does satisfy Condition $\text{PE}(V, c, \alpha, b, C)$ with $V(x) = (x+1)^{2\lambda}$ and $\alpha = (\lambda - 1)/2\lambda$. Since the upward jumps of $V(X)$ are clearly bounded for $X \geq 1$, the chain is tame by Theorem 21.

4.3. *Random walk Metropolis–Hastings.* For a more practical example, consider a random walk Metropolis–Hastings algorithm on $\mathbb{R}^d$, with proposal density $q$ and target density $p$. Fort and Moulines [6] consider the case where $q$ is symmetric and compactly supported and $\log p(z) \sim -|z|^s$, $0 < s < 1$, as $|z| \to \infty$. (When $d = 1$, this class of target densities includes distributions with tails typically heavier than the exponential, such as the Weibull distributions; see [6] for more details.) They show that under these conditions, the Metropolis–Hastings algorithm converges at *any* polynomial rate. In particular, it is possible to choose a scale function $V$ such that the chain satisfies Condition PE with $\alpha > 3/4$. Therefore, by Theorem 22, this chain is tame.

4.4. *Random walk on a half-line.* For our final example of a tame chain, we consider Example 5.1 of Tuominen and Tweedie [23]. This is the random walk on $[0, \infty)$ given by

$$(35) \qquad X_{n+1} = (X_n + Z_{n+1})^+,$$

where $\{Z_n\}$ is a sequence of i.i.d. real-valued random variables. We suppose that $\mathbb{E}[Z] = -\mu < 0$ (so that 0 is a positive-recurrent atom) and that $\mathbb{E}[(Z^+)^m] = \mu_m < \infty$ for some integer $m \geq 2$.

We also assume that $\mathbb{E}[r^{Z^+}] = \infty$ for all $r > 1$, and claim that this forces $X$ to be subgeometrically ergodic. To see this, consider the chain $\hat{X}$ which uses the same downward jumps as $X$ but which stays still when $X$ increases. That is,

$$\hat{X}_{n+1} = (\hat{X}_n - Z^-_{n+1})^+.$$

Let $\tau_0$ be the first time that $X$ hits 0 and let $\hat{\tau}_0$ be the corresponding hitting time for $\hat{X}$. Note that for all $n > 0$,

$$(36) \qquad \mathbb{E}_x[\hat{X}_{n \wedge \hat{\tau}_0}] \geq x - \mathbb{E}_x[n \wedge \hat{\tau}_0]\hat{\mu},$$

where $\hat{\mu} := -\mathbb{E}[Z; Z \leq 0] > 0$. Now, the left-hand side of (36) is dominated by $x$, and $\mathbb{E}_x[\hat{\tau}_0] < \infty$, so letting $n \to \infty$ yields

$$(37) \qquad \mathbb{E}_x[\tau_0] \geq \mathbb{E}_x[\hat{\tau}_0] \geq x/\hat{\mu}.$$

Thus, for $r > 1$,

$$\begin{aligned} \mathbb{E}_0[r^{\tau_0}] &= r\mathbb{E}_0[\mathbb{E}_{X_1}[r^{\tau_0}]] \\ &\geq r\mathbb{E}_0[r^{\mathbb{E}_{X_1}[\tau_0]}] \\ &\geq r\mathbb{E}_0[r^{X_1/\hat{\mu}}] = \infty \qquad \text{by assumption.} \end{aligned}$$

Therefore, by Theorem 3, $X$ is not geometrically ergodic.

Now, [11] show that if $m \geq 2$ is an integer, then $X$ satisfies Condition PE with $V(x) = (x+1)^m$ and $\alpha = (m-1)/m$. Clearly, the upward jumps of $V(X)$ when $X \notin C$ are not necessarily bounded, so Theorem 21 cannot be applied. Furthermore, if $m \leq 4$, then $\alpha \leq 3/4$, so Theorem 22 cannot be applied. However, we now show that $X$ is still tame when $m = 2$ (and thus tame for all $m \geq 2$).

(i) First, assume that the law of $Z$ is concentrated on $[-z_0, \infty)$ for some $z_0 > 0$, so $\mu_2 = \mathbb{E}[(Z^+)^2] < \infty$. Then if $x \geq z_0$,

$$\mathbb{E}_x[(X_1+1)^2] = \mathbb{E}[(x+1+Z)^2]$$
$$= (x+1)^2 + 2(x+1)\mathbb{E}[Z] + \mathbb{E}[Z^2]$$
$$\leq (x+1)^2 - 2\mu(x+1) + (\mu_2 + z_0^2).$$

Thus for any $0 < \beta < 1$, there exist $z_\beta > z_0$ and $b_\beta < \infty$ such that, with $V(x) = (x+1)^2$ and $\alpha = 1/2$,

(38) $$\mathbb{E}_x[V(X_1)] \leq V(x) - (2-\beta)\mu V^\alpha(x) + b_\beta \mathbf{1}_{[x \leq z_\beta]}.$$

Assume that $\beta < 1/4$ and a corresponding $z_\beta > z_0$ are fixed. Write $C_\beta = [0, z_\beta]$ and for $V(x) > z_\beta$, define $F(V(x)) = \lceil V^{1/2}(x)/\mu \rceil$. Iterating the drift condition (38), we obtain for $x \notin C$, with $F_x = F(V(x))$,

$$\mathbb{E}_x[V(X_{F_x})] \leq V(x) - (2-\beta)\mu \sum_{k=0}^{F_x-1} \mathbb{E}_x[V^{1/2}(X_k)] + b_\beta F_x$$

$$\leq (x+1)^2 - (2-\beta)\mu \sum_{k=0}^{F_x-1}(x+1-k\mu) + b_\beta F_x$$

(39) $$\leq \left(1 - (2-\beta) + \frac{(2-\beta)}{2}\right)(x+1)^2 + \gamma x \quad \text{since } \mathbb{E}_x[V^{1/2}(X_k)] = \mathbb{E}_x[(X_k+1)] \geq x+1-k\mu,$$

for some $\gamma > 0$,

$$\leq \frac{\beta}{2}V(x) + \gamma V^{1/2}(x).$$

Thus there exists a sublevel set $C'$ and a constant $b' < \infty$ such that if

$$F(x) = \begin{cases} \lceil x^{1/2}/\mu \rceil, & x \notin C', \\ 1, & x \in C', \end{cases}$$

then we obtain

$$\mathbb{E}_x[V(X_{F_x})] \leq \beta V(x) + b'\mathbf{1}_{C'}(x)$$

with $\beta < 1/4$. Since $\alpha = 1/2$, we satisfy $\log \beta < (1-\alpha)^{-1} \log \alpha$, so this chain is indeed tame.

(ii) In the general case, we can proceed by truncating the law of $Z$ at a level $-z_0$ so that the truncated distribution has a negative mean. The resulting chain, $X^*$ say, is tame by the above argument. However, $X^*$ stochastically dominates $X$ on the whole of $[0,\infty)$, so $X$ must also be tame.

A polynomial drift condition can still be shown to hold when $m \in (1, 2)$ [corresponding to drift $\alpha \in (0, 1/2)$]. Furthermore, it is quite simple to produce an adaptive subsampling scheme in this situation that produces a chain satisfying condition $GE(V, \beta, b, C)$. However, it is also necessary to make $\beta$ sufficiently small to satisfy part (b) of Definition 14 and we have not yet been able to achieve this. Therefore, it is unclear at present whether such chains are in fact tame.

**5. Conclusions and questions** We have introduced the concept of a *tame* Markov chain and shown that a *domCFTP* algorithm exists for all such chains. This algorithm is not expected to be practical in general, but it directly extends the results of [7] and [14]. In a practical setting, of course, one would use a dominating process that is better suited to the chain of interest. We have proven two sufficient conditions for a polynomially ergodic chain to be tame and provided an example which demonstrates that neither of these sufficient conditions are necessary.

Our suspicion, which is shared by those experts with whom we have discussed this, is that the following conjecture is true.

CONJECTURE 23. *There exists a chain satisfying Condition PE which is wild.*

On the other hand, we do not rule out the possibility that all polynomially ergodic chains are tame. A resolution of this conjecture would do much to further our understanding of such chains. The tame/wild classification provides some structure to the class of subgeometrically ergodic Markov chains that goes beyond the rate at which they converge to equilibrium. Although purely theoretical at present, this may prove to be important in understanding elaborate MCMC implementations: for a tame chain, the existence of a time change which produces a geometrically ergodic chain could possibly be exploited to improve the behavior of an MCMC algorithm.

It is also natural to ask what can be said about the more general case of subgeometric ergodicity. The drift condition

$$(40) \qquad PV \leq V - \phi \circ V + b\mathbf{1}_C$$

[where $\phi > 0$ is a concave, nondecreasing, differentiable function with $\phi'(t) \to 0$ as $t \to \infty$] is a generalization of (7) which can be shown to imply subgeometric ergodicity [5]. Much of the work in this paper extends naturally to chains satisfying this drift condition (see [2] for details). However, it is possible to produce a version $X$ of the Epoch chain of Section 4.1 that satisfies (40) but not (7). Furthermore, no subsampling scheme defined using a function $F$ of the form (16) will result in a geometrically ergodic chain, so this $X$ is wild. The existence of a perfect simulation algorithm for this and similar chains is also an open question.

DEPARTMENT OF STATISTICS
UNIVERSITY OF WARWICK
COVENTRY CV4 7AL
UNITED KINGDOM
E-MAIL: s.b.connor@warwick.ac.uk
    w.s.kendall@warwick.ac.uk